\documentclass[letterpaper,11pt]{article}
\usepackage{amsmath,amssymb,amsxtra,amsthm}
\usepackage{graphicx}

\theoremstyle{plain}

\theoremstyle{definition}

\newcommand{\gd}{\delta}
\newcommand{\eps}{\varepsilon}

\newcommand{\gq}{\vartheta}

\newcommand{\gz}{\zeta}

\newcommand{\tbf}{\textbf}

\title{New Proof of the Equation $$\sum_{k=1}^\infty \frac{\mu(k)}{k}=0.$$}
\author{Edmund Landau\\ Translated by Michael J.~Coons\\ Preliminary Version: September 2007}
\date{}
\begin{document}
\maketitle

\begin{center}{\large
[Neuer Beweis der Gleichung $\sum \mu(k)/k=0$,\\ {\em Inaugural-Dissertation}, Berlin, 1899.]}\\

Translation \copyright\ M.~J.~Coons 2007.
\end{center}

\newpage

\begin{center}

\tbf{{\Huge New Proof of the Equation $$\sum_{k=1}^\infty \frac{\mu(k)}{k}=0.$$
}}
\noindent\rule{5cm}{.5pt}

\vspace{1cm}

{\large \tbf{Inaugural--Dissertation}}\\
\vspace{.3cm}
for the \\
\vspace{.3cm}
acquisition of the title of Doctor from the Faculty of Philosophy\\
\vspace{.3cm}
of\\
\vspace{.3cm}
{\large\tbf{Friedrich--Wilhelms University in Berlin}}\\
\vspace{.3cm}
defended publicly and approved together with the attached theses\\
\vspace{.3cm}
on 15 July, 1899\\
\vspace{.3cm}
by\\
\vspace{.3cm}
{\Large \tbf{Edmund Landau}}

of Berlin

\vspace{1cm}
Opponents:
\end{center}

\tbf{Mr. Rudolf Zeigel}, Student of Mathematics

\tbf{Mr. Fritz Hartoge}, Student of Mathematics

\tbf{Ernst Steinitz}, Ph.D. , Privatdozent at the Royal Technical high school in Charlottenburg.

\vspace{1cm}

\begin{center}
\rule{7cm}{1pt}

\vspace{.5cm}
{\Large \tbf{Berlin 1899.}}
\end{center}

\newpage

\begin{center}
$$ $$
$$ $$
$$ $$
$$ $$
\tbf{{\Huge For my dear parents.}}
\end{center}

\newpage

$$ $$
$$ $$
$$ $$
$$ $$

The function $\mu(k)$ is normally defined as the number--theoretic function for which 
\begin{enumerate}
\item $\mu(1)=1$,
\item $\mu(k)=0$ when $k>1$ is divisible by a square,
\item $\mu(k)=(-1)^r$ when $k$ is the product of $r$ distinct primes.
\end{enumerate}

This statement was first expressed by Euler, that $$\sum_{k=1}^\infty \frac{\mu(k)}{k}=0$$ holds; that is, $\lim_{x\to\infty}\sum_{k=1}^x\frac{\mu(k)}{k}$ exists and equals $0$, the recent proof of which is due to von Mangoldt\footnote{) ``Beweis der Gleichung $\sum_{k=1}^\infty\frac{\mu(k)}{k}=0$''; Proceedings of the Royal Prussian Academy of Science of Berlin, 1897, pp. 835--852.}). The same goes for the investigations of Hadamard and de la Vall\'ee Poussin over the Riemann $\gz$--function, and it seems also, that without the use of these works, the present means of analysis is not enough to give a proof of Euler's statement. However, if one expects the results of those investigations to be in agreement with those of von Mangoldt, then one, as will be executed in the following, can arrive at the target along a quite short path.

The proof, which forms the contents of this dissertation, uses first the theorem\footnote{) This theorem is proven without the use of von Mangoldt's proof.}) of Hadamard\footnote{) Bulletin de la soci\'et\'e math\'ematique de France, Volume 24, 1896, p. 217.}) and de la Vall\'ee Poussin\footnote{) Annales de la soci\'et\'e scientifique de Bruxelles, Volume 20, Part 2, p. 251.}): 

``If $\gq(x):=\sum_{p\leq x}\log p$, then $\lim_{x\to\infty} \frac{\gq(x)}{x}=1.$'' However, apart from the use of this theorem, it is as elementary as can be for such a transcendent statement.

\section{}

Denote by $g(x)$ the sum $\sum_{k=1}^{[x]}\frac{\mu(k)}{k}$, where $[x]$ denotes the greatest integer less than or equal to $x$; more simply we write \begin{equation}g(x)=\sum_{k=1}^{x}\frac{\mu(k)}{k},\end{equation} where $k$ ranges over all positive integers less than or equal to $x$. The sum has meaning only for $x\geq 1$; thus, for $x<1$, set $g(x)=0$.

With the above notation, we read the two lemmas, which von Mangoldt proves in a simple way at the start of his paper\footnote{) 1. c., pp. 837--839.}) and which are also applied in the following one, as:

For all $x$ \begin{equation}|g(x)|\leq 1\ \footnote{) This lemma was, as indicated there, proved in the writings of Gram: ``Unders\o gelser abgaaende Maengden af Primtal under en given Graense,'' K.~Danske Videnskabernes Selskabs Skrifter, 6te Raekke, naturvidenskabelig og mathematisk Afdeling, II, 1884, pp. 197--198.})\end{equation} and for all $x\geq 1$ \begin{equation} \left| \log x\cdot g(x)-\sum_{k=1}^x \frac{\mu(k)\log k}{k}\right| \leq 3+\gamma\end{equation} where $\gamma$ denotes Euler's constant.

The inequality (3), which von Mangoldt only derived in order to apply it in a certain place in his proof\footnote{) 1. c., p. 843.}), serves in the following one as the basis of the whole investigation.

Concerning the sum $\sum_{k=1}^x \frac{\mu(k)\log k}{k}$, M\"obuis\footnote{) ``\"Uber eine besondere Art von Umkehrung Reihen,'' Journal f\"ur die reine und angewandte Mathematik, Volume 9, p. 122.}) believed he had proved that for sufficiently large $x$, its difference from $-1$ is arbitrarily small: however, his proof is not sound. Although new writers consider it probable\footnote{) E.g., Mertens proved, which the general validity of the inequality condition assumes: $$\left|\sum_{k=1}^x \mu(k)\right|\leq \sqrt{x},$$ that this theorem is correct, if this relation is generally fulfilled (Proceedings of the Vienna Academy, math.-nat. Kl., Volume 106, Dept. 2a, p. 774.)}) that $$\lim_{x=\infty}\sum_{k=1}^x\frac{\mu(k)\log k}{k}$$ exists and equals $-1$, it has yet to be proven that for all $x$, $\sum_{k=1}^x\frac{\mu(k)\log k}{k}$ is contained between two finite boundaries. Now since (3) yields $$\left| g(x)-\frac{1}{\log x}\sum_{k=1}^x \frac{\mu(k)\log k}{k}\right| \leq \frac{3+\gamma}{\log x},$$ it follows, with use of the Euler--v.~Mangoldt Theorem, that $$\lim_{x=\infty} g(x)=0$$ so that $$\frac{1}{\log x}\sum_{k=1}^x \frac{\mu(k)\log k}{k}$$ approaches $0$ as $x\to\infty$.

If, in reverse, it was successfully proven that $$\lim_{x=\infty}\frac{1}{\log x}\sum_{k=1}^x \frac{\mu(k)\log k}{k}$$ exists and equals 0, then one would trivially have that $$\sum_{k=1}^x \frac{\mu(k)}{k}=0,$$ since for any $\gd$ there is a $G$ such that for all $x\geq G$ $$\left|\frac{1}{\log x}\sum_{k=1}^x \frac{\mu(k)\log k}{k}\right|\leq \frac{\gd}{2}$$ and $$0<\frac{3+\gamma}{\log x}\leq \frac{\gd}{2},$$ thus it follows for $x\geq G$: 
\begin{align*} |g(x)| &= \left|\left(g(x)-\frac{1}{\log x}\sum_{k=1}^x \frac{\mu(k)\log k}{k}\right)+\frac{1}{\log x}\sum_{k=1}^x \frac{\mu(k)\log k}{k}\right|\\
&\leq \left|\left(g(x)-\frac{1}{\log x}\sum_{k=1}^x \frac{\mu(k)\log k}{k}\right)\right|+\left|\frac{1}{\log x}\sum_{k=1}^x \frac{\mu(k)\log k}{k}\right|\\
&\leq \frac{3+\gamma}{\log x}+\left|\frac{1}{\log x}\sum_{k=1}^x \frac{\mu(k)\log k}{k}\right|\\
&\leq \frac{\gd}{2}+\frac{\gd}{2}=\gd,\end{align*} also $$\lim_{x=\infty} g(x)=0.$$

The proof, that for \begin{equation} f(x)=\sum_{k=1}^x\frac{\mu(k)\log k}{k},\end{equation} one has \begin{equation} \lim_{x=\infty}\frac{f(x)}{\log x}=0\end{equation} will be supplied in what follows.

In order not to have to interrupt the course of the investigation, we note the following simple lemma, which was already known to Gram\footnote{) 1. c., p. 197, where separately for all valid $r$ in equation (43) set $r=1$.}): it is \begin{equation} \sum_{\nu=1}^x \frac{1}{\nu}g\left(\frac{x}{\nu}\right) = g(x)+\frac{1}{2}g\left(\frac{x}{2}\right)+\frac{1}{3}g\left(\frac{x}{3}\right)+\cdots+\frac{1}{[x]}g\left(\frac{x}{[x]}\right)=1.\end{equation} The $\nu$--th summand $\frac{1}{\nu}g\left(\frac{x}{\nu}\right)$ contains the sum of the terms $$\frac{1}{\nu}\frac{\mu(1)}{1}=\frac{\mu(1)}{\nu},\ \frac{1}{\nu}\frac{\mu(2)}{2}=\frac{\mu(2)}{2\nu},\cdots,\ \frac{1}{\nu}\frac{\mu(n)}{n}=\frac{\mu(n)}{n\nu},\cdots,\ \frac{1}{\nu}\frac{\mu\left[\frac{x}{\nu}\right]}{\left[\frac{x}{\nu}\right]}=\frac{\mu\left[\frac{x}{\nu}\right]}{\left[\frac{x}{\nu}\right]\nu};$$ the sum $\sum_{\nu=1}^x\frac{1}{\nu}g\left(\frac{x}{\nu}\right)$ consists of terms of the form $\frac{\mu(n)}{t}$, where $n$ is a divisor of $t$, and $t$ runs through the integers from 1 to $[x]$; that is, $$\sum_{\nu=1}^x\frac{1}{\nu}g\left(\frac{x}{\nu}\right)=\sum_{t=1}^x\frac{1}{t}\sum_{n|t}\mu(n);$$ now since $\sum_{n|t}\mu(n)$ is 1 for $t=1$ and $0$ otherwise, we have $$\sum_{\nu=1}^x\frac{1}{\nu}g\left(\frac{x}{\nu}\right)=1.$$

\section{}

If one lets $x=p_1p_2\cdots p_r$ in the defining equation (4) of $f(x)$, replaces\footnote{) The $\log k$ is multiplied by the factor $\frac{\mu(k)}{k}$, which occurs only for $k$ that are the product of distinct primes.}) $\log(x)$ by $\log p_1+\cdots+\log p_r$ and gathers like terms in which the logarithm is applied to the same prime number, then (4) becomes an equation of the form $$f(x)=\sum F(p,x) \log p$$ where the sum extends over all prime numbers $p\leq x$. As was easily given by von Mangoldt\footnote{) 1. c., p. 840.}) for another purpose, $$F(p,x)=-\left(\frac{1}{p}\sum_{k=1}^{\frac{x}{p}}\frac{\mu(k)}{k}+\frac{1}{p^2}\sum_{k=1}^{\frac{x}{p^2}}\frac{\mu(k)}{k}+\frac{1}{p^3}\sum_{k=1}^{\frac{x}{p^3}}\frac{\mu(k)}{k}+\cdots\right),$$ a series which has only a finite number of non-zero summands, since the summation index of the $i$--th sum runs from 1 to $\left[\frac{x}{p^i}\right]$, so that $p^i\leq x$ gives $i\leq\frac{\log x}{\log p}$. Therefore $$f(x)=-\sum_{p\leq x}\log p\left(\frac{1}{p}g\left(\frac{x}{p}\right)+\frac{1}{p^2}g\left(\frac{x}{p^2}\right)+\frac{1}{p^3}g\left(\frac{x}{p^3}\right)+\cdots\right),$$ \begin{equation} f(x)=-\sum_{p\leq x}\frac{\log p}{p}g\left(\frac{x}{p}\right)-\sum_{p\leq x}\log p\left(\frac{1}{p^2}g\left(\frac{x}{p^2}\right)+\frac{1}{p^3}g\left(\frac{x}{p^3}\right)+\cdots\right).\end{equation}

According to (2), for all $y$ $$|g(y)|\leq 1,$$ so that the absolute value of the second sum in (7) is \begin{align*} &\left|\sum_{p\leq x}\log p\left(\frac{1}{p^2}g\left(\frac{x}{p^2}\right)+\frac{1}{p^3}g\left(\frac{x}{p^3}\right)+\cdots\right)\right|\\
&\leq \sum_{p\leq x}\log p\left(\frac{1}{p^2}\left|g\left(\frac{x}{p^2}\right)\right|+\frac{1}{p^3}\left|g\left(\frac{x}{p^3}\right)\right|+\cdots\right)\\
&\leq \sum_{p\leq x}\log p\left(\frac{1}{p^2}+\frac{1}{p^3}+\frac{1}{p^4}+\cdots+\frac{1}{p^n}+\cdots\right)\\
&\leq \sum_{p\leq x}\log p\left(\frac{1}{p^2}+\frac{1}{2p^2}+\frac{1}{2^2p^2}+\cdots+\frac{1}{2^{n-2}p^2}+\cdots\right)\\
&\leq \sum_{p\leq x}\frac{\log p}{p^2}\left(1+\frac{1}{2}+\frac{1}{4}+\frac{1}{8}+\cdots\right)\\
&\leq 2\sum_{p\leq x}\frac{\log p}{p^2}\\
&<2\sum_{\nu=1}^\infty\frac{\log\nu}{\nu^2}.\end{align*} It is well known that $\sum_{\nu=1}^\infty\frac{\log\nu}{\nu^2}$ is convergent; thus as $x\to\infty$ the sum on the right-hand side of (7) approaches either a certain bound, or its value oscillates between two finite uncertain bounds. In either case, the quotient with $\log x$ approaches $0$ as $x\to\infty$.

Denote by $h(x)$, the function defined by \begin{equation} h(x)=\sum_{p\leq x}\frac{\log p}{p}g\left(\frac{x}{p}\right).\end{equation} If $\lim_{x=\infty}\frac{h(x)}{\log x}$ exists and equals $0$, then according to (7) $$\frac{f(x)+h(x)}{\log x}=\frac{-1}{\log x}\sum_{p\leq x}\log p\left(\frac{1}{p^2}g\left(\frac{x}{p^2}\right)+\frac{1}{p^3}g\left(\frac{x}{p^3}\right)+\cdots\right),$$ and as we just saw, as $x$ gets large the right-hand side approaches $0$, from which the correctness of statement (5) follows.

The proof of Euler's statement, that $$\sum_{k=1}^x\frac{\mu(k)}{k}=0,$$ thus depends on the proof of the statement $$\lim_{x=\infty}\frac{h(x)}{\log x}=0,$$ which will be furnished in the following section.

\section{}

Recall that the function $\gq(x)$ \footnote{) See the theorem of Hadamard and de la Vall\'ee Poussin in the introduction.}) is defined for all positive $\nu$ by $$\gq(\nu)-\gq(\nu-1) =\begin{cases} \log \nu &\mbox{ if $\nu$ is prime,}\\
0 &\mbox{ if $\nu$ is composite or 1},\\
\log \nu=0 &\mbox{ if $\nu=1$.}\end{cases}$$

And $$h(x)=\sum_{\nu=1}^x \frac{\gq(\nu)-\gq(\nu-1)}{\nu} g\left(\frac{x}{\nu}\right),$$ where the sum ranges over the integers between $1$ and $x$.

In the place of the function $\gq(x)$, use \begin{equation}\gq(x)=x\{1+\eps(x)\},\end{equation} where the function $\eps(x)$ takes only non-negative values of $x$, and $\eps(0)=0$. We note the following properties of $\eps(x)$:

\begin{enumerate} 
\item Since by definition, $\gq(x)$ is never negative, then always $$\eps(x)\geq -1.$$
\item As shown by Mertens\footnote{) ``Ein Beitrag zur analytischen Zahlentheorie,'' Journal f\"ur die reine und angewandte Mathematik, Volume 78, p. 48.}), for all $x$ $$\gq(x)<2x,$$ so that always $$\eps(x)<1;$$ therefore, we gain the inequality \begin{equation} |\eps(x)|\leq 1.\end{equation}
\item The theorem cited in the introduction\footnote{) See the theorem of Hadamard and de la Vall\'ee Poussin in the introduction.}), that $$\lim_{x=\infty}\frac{\gq(x)}{x}=1,$$ gives \begin{equation} \lim_{x=\infty}\eps(x)=0.\end{equation}
\end{enumerate}

The introduction of the function $\eps(x)$ yields, for $h(x)$, \begin{align*} h(x)&= \sum_{\nu=1}^x \frac{\nu+\nu\eps(\nu)-(\nu-1)-(\nu-1)\eps(\nu-1)}{\nu}g\left(\frac{x}{\nu}\right)\\
&= \sum_{\nu=1}^x \left\{\frac{1}{\nu}g\left(\frac{x}{\nu}\right)+\left(\eps(\nu)-\frac{\nu-1}{v}\eps(\nu-1)\right)g\left(\frac{x}{\nu}\right)\right\}\\
&= \sum_{\nu=1}^x \frac{1}{\nu}g\left(\frac{x}{\nu}\right)+\sum_{\nu=1}^x \left(\eps(\nu)-\eps(\nu-1)+\frac{1}{\nu}\eps(\nu-1)\right)g\left(\frac{x}{\nu}\right).\end{align*}

Using Eqs.~(8) and (9), $$\sum_{\nu=1}^x \frac{1}{\nu}g\left(\frac{x}{\nu}\right)=1;$$ yielding \begin{equation} h(x)-1=\sum_{\nu=1}^x \left(\eps(\nu)-\eps(\nu-1)\right)g\left(\frac{x}{\nu}\right)+\sum_{\nu=1}^x \frac{1}{\nu}\eps(\nu-1)g\left(\frac{x}{\nu}\right).\end{equation}

For the first of the two sums in (12) we get \begin{align*} \sum_{\nu=1}^x (\eps(\nu)&-\eps(\nu-1))g\left(\frac{x}{\nu}\right)\\ &= \sum_{\nu=1}^x \eps(x)\left(g\left(\frac{x}{\nu}\right)-g\left(\frac{x}{\nu+1}\right)\right)+\eps([x])g\left(\frac{x}{[x]+1}\right)\\
&=\sum_{\nu=1}^x \eps(\nu)\left(g\left(\frac{x}{\nu}\right)-g\left(\frac{x}{\nu+1}\right)\right),\end{align*} where $x<[x]+1$ so that $g\left(\frac{x}{[x]+1}\right)=0$.

If in the second sum in (12) we write $\nu+1$ in place of $\nu$, we have $$
\sum_{\nu=0}^{x-1} \frac{1}{\nu+1}\eps(\nu)g\left(\frac{x}{\nu+1}\right)=\sum_{\nu=1}^{x-1} \frac{1}{\nu+1}\eps(\nu)g\left(\frac{x}{\nu+1}\right)$$ and so \begin{equation} h(x)-1=\sum_{\nu=1}^x \eps(\nu)\left(g\left(\frac{x}{\nu}\right)-g\left(\frac{x}{\nu+1}\right)\right)+\sum_{\nu=1}^{x-1} \frac{1}{\nu+1}\eps(\nu)g\left(\frac{x}{\nu+1}\right).\end{equation}

Let $\gd$ be an arbitrary small positive quantity. Then by (11), there is a $G$ such that for all $\nu\geq G$ \begin{equation} |\eps(\nu)|\leq \frac{\gd}{3}.\end{equation} This yields $$\left|\sum_{\nu=1}^x \eps(\nu)\left(g\left(\frac{x}{\nu}\right)-g\left(\frac{x}{\nu+1}\right)\right)\right|\leq \left|\sum_{\nu=1}^{G-1} \eps(\nu)\left(g\left(\frac{x}{\nu}\right)-g\left(\frac{x}{\nu+1}\right)\right)\right|$$
$$+\left|\sum_{\nu=G}^x \eps(\nu)\left(g\left(\frac{x}{\nu}\right)-g\left(\frac{x}{\nu+1}\right)\right)\right|$$
$$\leq \sum_{\nu=1}^{G-1} |\eps(\nu)|\left|\left(g\left(\frac{x}{\nu}\right)-g\left(\frac{x}{\nu+1}\right)\right)\right|+\sum_{\nu=G}^x |\eps(\nu)|\left|\left(g\left(\frac{x}{\nu}\right)-g\left(\frac{x}{\nu+1}\right)\right)\right|.$$ As the right-hand side is increased, $|\eps(\nu)|$ goes to 1 in the first sum (by (10)), and goes to $\gd/3$ in the second sum (by (14)), yielding \begin{equation} 
\left|\sum_{\nu=1}^x \eps(\nu)\left(g\left(\frac{x}{\nu}\right)-g\left(\frac{x}{\nu+1}\right)\right)\right|\leq \sum_{\nu=1}^{G-1} \left|\left(g\left(\frac{x}{\nu}\right)-g\left(\frac{x}{\nu+1}\right)\right)\right|\end{equation} $$+\frac{\gd}{3}\sum_{\nu=G}^{x} \left|\left(g\left(\frac{x}{\nu}\right)-g\left(\frac{x}{\nu+1}\right)\right)\right|.$$ Now \begin{align*} \sum_{\nu=1}^{G-1} \left|\left(g\left(\frac{x}{\nu}\right)-g\left(\frac{x}{\nu+1}\right)\right)\right| &\leq \sum_{\nu=1}^{G-1} \left\{\left|\left(g\left(\frac{x}{\nu}\right)\right|-\left|g\left(\frac{x}{\nu+1}\right)\right)\right|\right\}\\
&= \sum_{\nu=1}^{G-1}(1+1)\ \mbox{ (by (2))},\end{align*} so that \begin{equation} \sum_{\nu=1}^{G-1} \left|\left(g\left(\frac{x}{\nu}\right)-g\left(\frac{x}{\nu+1}\right)\right)\right|\leq 2(G-1)\end{equation} and \begin{align*} \sum_{\nu=G}^{x} \left|\left(g\left(\frac{x}{\nu}\right)-g\left(\frac{x}{\nu+1}\right)\right)\right| &\leq \sum_{\nu=1}^{x} \left|\left(g\left(\frac{x}{\nu}\right)-g\left(\frac{x}{\nu+1}\right)\right)\right|\\
&= \sum_{\nu=1}^{x}\left| \sum_{k=1}^{\frac{x}{\nu}}\frac{\mu(k)}{k}-\sum_{k=1}^{\frac{x}{\nu+1}}\frac{\mu(k)}{k}\right| \\
&= \sum_{\nu=1}^{x}\left| \sum \frac{\mu(k)}{k}\right|,\end{align*} where $k$ ranges over all integers in the interval $\left(\frac{x}{\nu+1},\frac{x}{\nu}\right]$. Therefore 
$$\sum_{\nu=G}^{x} \left|\left(g\left(\frac{x}{\nu}\right)-g\left(\frac{x}{\nu+1}\right)\right)\right|\leq \sum_{\nu=1}^{x}\sum_{\frac{x}{\nu}\geq k>\frac{x}{\nu+1}}\frac{|\mu(k)|}{k}\leq \sum_{\nu=1}^{x}\sum_{\frac{x}{\nu}\geq k>\frac{x}{\nu+1}}\frac{1}{k}$$ $$=\sum_{x\geq k>\frac{x}{2}}\frac{1}{k}+\sum_{\frac{x}{2}\geq k>\frac{x}{3}}\frac{1}{k}+\sum_{\frac{x}{3}\geq k>\frac{x}{4}}\frac{1}{k}+\cdots+\sum_{\frac{x}{[x]-1}\geq k>\frac{x}{[x]}}\frac{1}{k}$$ $$+\sum_{\frac{x}{[x]}\geq k>1}\frac{1}{k}=\sum_{k=1}^x\frac{1}{k},$$ and since always $$\sum_{k=1}^x\frac{1}{k}\leq \log x+1,$$ we have \begin{equation}\sum_{\nu=G}^{x} \left|\left(g\left(\frac{x}{\nu}\right)-g\left(\frac{x}{\nu+1}\right)\right)\right|\leq \log x+1. \end{equation}

Replacing both sums of the right-hand side of inequality (15) by the results gained in (16) and (17) yields \begin{equation} \left|\sum_{\nu=1}^{x} \eps(x)\left(g\left(\frac{x}{\nu}\right)-g\left(\frac{x}{\nu+1}\right)\right)\right|\leq 2(G-1)+\frac{\gd}{3}(\log x+1).\end{equation}

The handling of the second sum in (13) is somewhat simpler. We have \begin{align*} \left|\sum_{\nu=1}^{x-1} \frac{1}{\nu+1}\eps(\nu)g\left(\frac{x}{\nu+1}\right)\right| &\leq \sum_{\nu=1}^{x-1} \frac{1}{\nu+1}|\eps(\nu)|\left|g\left(\frac{x}{\nu+1}\right)\right|\\
&\leq \sum_{\nu=1}^{x-1} \frac{|\eps(\nu)|}{\nu+1}\\
&= \sum_{\nu=1}^{G-1} \frac{|\eps(\nu)|}{\nu+1}+\sum_{\nu=G}^{x-1} \frac{|\eps(\nu)|}{\nu+1}\\
&\leq \sum_{\nu=1}^{G-1} \frac{1}{\nu+1}+\sum_{\nu=G}^{x-1} \frac{1}{\nu+1}\frac{\gd}{3}\\
&\leq  \sum_{\nu=1}^{G-1} 1+\frac{\gd}{3}\sum_{\nu=1}^{x-1} \frac{1}{\nu},\end{align*} so that \begin{equation} \left|\sum_{\nu=1}^{x-1} \frac{1}{\nu+1}\eps(\nu)g\left(\frac{x}{\nu+1}\right)\right|\leq G-1+\frac{\gd}{3}(\log x+1).\end{equation}

With help from the inequalities (18) and (19), (13) becomes $$|h(x)|\leq 1+2G-2 +\frac{\gd}{3}\log x+\frac{\gd}{3}+G-1+\frac{\gd}{3}\log x+\frac{\gd}{3}$$ $$=3G-2+\frac{2}{3}\gd+\frac{2}{3}\gd\log x,$$ thus for all $$x\geq e^{\frac{3G-2+\frac{2}{3}\gd}{\frac{1}{2}\gd}},$$ we have $$3G-2+\frac{2}{3}\gd\leq \frac{1}{3}\gd\log x,$$ so that $$|h(x)|\leq \frac{1}{3}\gd\log x+\frac{2}{3}\gd\log x=\gd\log x,$$ which yields \begin{equation} \left|\frac{h(x)}{\log x}\right|\leq \gd.\end{equation}

For such a $\gd$ there is a $\xi$ assignable, such that for all $x\geq \xi,$ (20) is fulfilled; therefore the $\lim_{x=\infty}\frac{h(x)}{\log x}$ exists and equals 0. Thus all the results shown in the first two paragraphs of this work are valid; that is, the $\lim_{x=\infty}\sum_{k=1}^x \frac{\mu(k)}{k}$ exists and equals 0,  and thus the correctness of the equation named in the title, briefly $$\sum_{k=1}^x \frac{\mu(k)}{k}=0.$$

\section{}

If we define\footnote{) von Mangoldt, 1. c., p. 850.}) $$M(x)=\sum_{k=1}^x \mu(k),$$ then with help of the proven result, $$\lim_{x=\infty} g(x)=0,$$ we have $$\lim_{x=\infty}\frac{M(x)}{x}=0.$$ Von Mangoldt\footnote{) 1. c., pp. 849--851.}) proved this indirectly by use of the identity $$g(x)=\sum_{k=1}^x \frac{\mu(k)}{k}=\sum_{k=1}^x (M(k)-M(k-1))\frac{1}{k}.$$ It can be furnished as follows directly. From the equation $$M(x)=\sum_{k=1}^x \mu(k)=\sum_{k=1}^x \frac{\mu(k)}{k}\cdot k=\sum_{k=1}^x (g(k)-g(k-1))\cdot k,$$ it follows that \begin{align*} M(x)&=\sum_{k=1}^{x-1} g(k)(k-(k+1))+g(x)[x]\\
&=-\sum_{k=1}^{x-1} g(k)+g(x)[x],\end{align*} so that since for $\gd>0$ there is a $G$ such that for all $k\geq G$ $$|g(k)|\leq \frac{\gd}{3}$$ for all $x\geq G$ \begin{align*} |M(x)|&\leq \sum_{k=1}^{G-1}|g(k)|+\sum_{k=G}^{x-1}|g(k)|+|g(x)|\cdot x\\
&\leq G-1+\frac{\gd}{3}([x]-G)+\frac{\gd}{3}x,\end{align*} so that $$\left|\frac{M(x)}{x}\right|\leq \frac{G-1-\frac{\gd}{3}G}{x}+\frac{2}{3}\gd,$$ then for $$x\geq \frac{G-1-\frac{\gd}{3}G}{\frac{\gd}{3}}$$ and at the same time greater than or equal to $G$, $$\left|\frac{M(x)}{x}\right|\leq \frac{1}{3}\gd+\frac{2}{3}\gd=\gd,$$ with which the statement $$\lim_{x=\infty}\frac{M(x)}{x}=0$$ is proved.\\

\begin{center}
\rule{7cm}{1pt}
\end{center}

\newpage

\section*{Theses}

\begin{enumerate}
\item It is desirable during every existence proof of a mathematical quantity to be led, at the same time on the way to the result, to the actual existing quantity.

\item A boundary between arithmetic and analytic areas of mathematics cannot be drawn.

\item The concept of the semiconvergent series is a relative concept.

\item Out of the impossibility of perpetual motion of second kind comes the proof of the second law of thermodynamics.

\item It did not succeed, the justifying of psychology on an exactly mathematical basis.
\end{enumerate}



\end{document}